\documentclass[12pt]{amsart}
\usepackage{amsmath,amsthm,amsfonts,amssymb,eucal}

\renewcommand {\a}{ \alpha }

\newcommand{\y}{\eta}

\newcommand{\vare}{\varepsilon}
\newcommand{\g}{\gamma}
\newcommand{\G}{\Gamma}

\newcommand{\vark}{\varkappa}
\newcommand{\varf}{\varphi}
\renewcommand{\d}{\delta}

\newcommand{\s}{\sigma}
\renewcommand{\l}{\lambda}

\newcommand{\p}{\partial}

\newcommand{\Om}{\Omega}

\newcommand{\R}{ \mathbb R}
\newcommand{\C}{ \mathbb C}
\newcommand{\N}{ \mathbb N}

\newcommand{\CL}{\mathcal L}

\newcommand{\CC}{\mathcal C}
\newcommand{\CD}{\mathcal D}
\newcommand{\CE}{\mathcal E}
\newcommand{\CF}{\mathcal F}
\newcommand{\CG}{\mathcal G}

\newcommand{\CK}{\mathcal K}

\newcommand {\GH}{\mathfrak H}

\newcommand {\gm}{\mathfrak m}

\newcommand {\ba}{\mathbf a}
\newcommand {\bd}{\mathbf d}

\newcommand {\bb}{\mathbf b}

\newcommand {\BA}{\mathbf A}
\newcommand {\BB}{\mathbf B}

\newcommand {\BI}{\mathbf I}

\newcommand {\BT}{\mathbf T}

\newcommand {\BX}{\mathbf X}

\newcommand{\CA}{\mathcal A}

\newcommand{\wh}{\widehat}

\DeclareMathOperator{\re}{Re} 
\DeclareMathOperator{\res}{\restriction}
\DeclareMathOperator{\Int}{Int} 
\DeclareMathOperator{\codim}{codim}

\newtheorem{thm}{Theorem}[section]

\newtheorem{lem}[thm]{Lemma}
\newtheorem{prop}[thm]{Proposition}

\theoremstyle{definition}
\newtheorem{defn}[thm]{Definition}

\theoremstyle{remark}

\numberwithin{equation}{section}

\newcommand{\thmref}[1]{Theorem~\ref{#1}}

\newcommand{\bsymb}{\boldsymbol}

\begin{document}

\title[spectrum of irreversible quantum graphs]
{Smilansky's model of irreversible quantum graphs, II: the point
spectrum}
\author[W.D. Evans]{W.D. Evans}
\address{School of Mathematics\\ Cardiff University\\
         23 Senghennydd Road\\ Cardiff CF24 4AG\\  UK}
\email{EvansWD@cardiff.ac.uk}
\author[M. Solomyak]{M. Solomyak}
\address{Department of Mathematics\\The Weizmann Institute of Science\\
Rehovot 76100\\Israel} \email{michail.solomyak@weizmann.ac.il}
\date{18th May, 2005}

\begin{abstract}
In the model suggested by Smilansky \cite{SM} one studies an
operator describing the interaction between a quantum graph and a
system of $K$ one-dimensional oscillators attached at different
points of the graph. This paper is a continuation of \cite{ESI} in
which we started an investigation of the case $K>1$. For the sake
of simplicity we consider $K=2$, but our argument applies to the
general situation. In this second part of the paper we apply the
variational approach to the study of the point spectrum.
\end{abstract}

\maketitle

\section{Introduction}\label{int}

In Smilansky's model of irreversible quantum graphs, the
interaction between a quantum graph and a finite system of
one-dimensional harmonic oscillators attached at various vertices
of the graph is studied. The paper \cite{SM} may be consulted for
the physical background and motivation, and \cite{Ku} for a survey
of recent work on quantum graphs. Our concern here is the spectral
analysis of the self-adjoint operator which generates the
dynamical system, and it suffices to have a precise description of
the analytic problem. This paper continues the study in \cite{ESI}
where a detailed description of the problem may be found and a
survey of earlier results in the literature given. As in
\cite{ESI}, we consider the case of two oscillators attached to
the graph constituted by $\R$ at vertices $\pm 1$. This special
case retains the main features of the general case without
obscuring the argument with technical complications.

 On a formal level, the problem is described by the differential
expression
\begin{equation}\label{e1.1}
    {\CA}U = -U_{x^2}'' + \frac{\nu_+^2}{2} (-U_{q_+^2}''+  q_+^2 U)
     + \frac{\nu_-^2}{2} (-U_{q_-^2}''+  q_-^2 U)
\end{equation}
for $x\in \R,\ q_{\pm} \in \R$, together with the following
`transmission', or `matching' conditions across the planes $x=\pm
1$ in ${\R}^3$:
\begin{equation}\label{e1.2}\begin{split}
   & U_x'(1+,q_+,q_-)-U_x'(1-,q_+,q_-)= \alpha_+ q_+U(0,q_+,q_-),\\
&U_x'(-1+,q_+,q_-)-U_x'(-1-,q_+,q_-)= \alpha_- q_-U(0,q_+,q_-).
\end{split}\end{equation}
The parameters $\alpha_{\pm}$ are real and can be assumed to be
non-negative since, for instance, replacing $\alpha_{+}$ by $-
\alpha_{+}$ corresponds to replacing $q_+$ by $-q_+$ and this has
no effect on the problem to be investigated. The parameters
$\nu_{\pm}$ are fixed positive numbers throughout. To shorten our
notation, we set $\bsymb{\a}=(\a_+,\a_-)$ and
$\bsymb{\nu}=(\nu_+,\nu_-)$.

Let $\chi_n, n\in {\N}_0$, be the normalized Hermite functions in
$L^2({\R}).$ The sequence $\{\chi_n\}_{ n\in {\N}_0}$ is then an
orthonormal basis in $L^2({\R})$ and any $U \in L^2({\R}^3)$ can
be written as
\begin{equation*}
    U(x,q_+,q_-)= \sum_{m,n \in {\N}_0}u_{m,n}(x)\chi_m (q_+)\chi_n(q_-)
\end{equation*}
for some $u_{m,n} \in L^2(\R)$. We write $U \sim \{u_{m,n}\}$ to
indicate this representation. The mapping $U\mapsto\{u_{m,n}\}$ is
an isometry of $\GH=L^2(\R^3)$ onto the Hilbert space
$\ell^2(\N_0^2;L^2(\R))$. For $U \sim \{u_{m,n}\}$ we have ${\CA}U
\sim \{L_{m,n}u_{m,n}\}$, where
\begin{equation}\label{sa.1}
(L_{m,n}u)(x)=-u''(x)+r_{m,n}u(x),\qquad x\neq\pm1;
\end{equation}
\begin{equation}\label{sa.2}
r_{m,n}=\nu_+^2(m+1/2)+\nu_-^2(n+1/2),\qquad m,n\in\N_0.
\end{equation}
The number
\[r_{0,0}=(\nu_+^2+\nu_-^2)/2\]
plays a special role since it appears in the formulations of all
our basic results.

 The conditions \eqref{e1.2} at $x=\pm1$ become
\begin{equation}\label{1}\begin{split}
\sum_{m,n \in {\N}_0}[u_{m,n}'](1)\chi_m(q_+) \chi_n(q_-) =
\sum_{m,n
\in {\N}_0}\a_+q_+ \chi_m(q_+) \chi_n(q_-), \\
\sum_{m,n \in {\N}_0}[u_{m,n}'](-1)\chi_m(q_+) \chi_n(q_-) =
\sum_{m,n \in {\N}_0}\a_-q_- \chi_m(q_+) \chi_n(q_-), \end{split}
\end{equation}
where we have used the notation
\[
[u'](a) := u'(a+0) - u'(a-0). \] On using the recurrence relation
\begin{equation*}
\sqrt{ k+1} \chi_{k+1}(q) - \sqrt 2 q\chi_k(q) +\sqrt k
\chi_{k-1}(q) =0, \ \ \ q \in \R,
\end{equation*}
the matching conditions (\ref{1}) reduce to
\begin{equation}\label{sa.3}\begin{split}
[u'_{m,n}](1)&=\frac{\a_+}{\sqrt
2}\left(\sqrt{m+1}u_{m+1,n}(1)+\sqrt{m}u_{m-1,n}(1)
\right);\\
[u'_{m,n}](-1)&=\frac{\a_-}{\sqrt
2}\left(\sqrt{n+1}u_{m,n+1}(-1)+\sqrt{n}u_{m,n-1}(-1)
\right).\end{split}
\end{equation}

The operator realization of (\ref{e1.1}) and (\ref{e1.2}) in the
Hilbert space $\GH$, which we denote by
$\BA_{\bsymb{\a},\bsymb{\nu}}$ can  now be defined. Its domain
$\CD_{\bsymb{\a},\bsymb{\nu}}$ is given by
\begin{defn}
 {\it An element $U\sim\{u_{m,n}\}$
lies in $\CD_{\bsymb{\a},\bsymb{\nu}}$ if and only if

\noindent  1. $u_{m,n}\in H^1(\R)$ for all $m,n$;

\noindent 2. for all $m,n$, the restriction of $u_{m,n}$ to each
interval $(-\infty,-1), (-1,1), (1,\infty),$ lies in $H^2$ and
moreover,
\begin{equation*}
\sum_{m,n}\int_{\R}|L_{m,n}u_{m,n}|^2dx<\infty;
\end{equation*}
 3. the conditions (\ref{sa.3}) are satisfied.}
\end{defn}

\vskip0.2cm Along with the set $\CD_{\bsymb{\a},\bsymb{\nu}}$, we
define its subset
\begin{equation*}
\CD_{\bsymb{\a},\bsymb{\nu}}^{\bullet}=\left\{
U\in\CD_{\bsymb{\a},\bsymb{\nu}}: U\sim\{u_{m,n}\} \
\rm{finite}\right\},
\end{equation*}
where by {\it finite} we mean that the sequence has only a finite
number of non-zero components. \vskip0.2cm

The operator $\BA_{\bsymb{\a},\bsymb{\nu}}$ in $\GH$ is defined on
the domain $\CD_{\bsymb{\a},\bsymb{\nu}}$ by
\begin{equation*}
    \BA_{\bsymb{\a},\bsymb{\nu}} U \sim \{L_{m,n}u_{m,n}\} \ \  {\rm{for}}\ \ U \sim
    \{u_{m,n}\} \in \CD_{\bsymb{\a},\bsymb{\nu}}
\end{equation*}
where $L_{m,n}$ is given by (\ref{sa.1}). We denote the
restriction of $\BA_{\bsymb{\a},\bsymb{\nu}}$ to
$\CD_{\bsymb{\a},\bsymb{\nu}}^{\bullet}$ by
$\BA_{\bsymb{\a},\bsymb{\nu}}^{\bullet}$.

The following statement is proved in \cite{ESI}, Theorem 2.3.
\begin{thm}\label{ES0}
The operator $\BA_{\bsymb{\a},\bsymb{\nu}}$ is self-adjoint for
all $\a_{\pm} \ge 0,$ and is the closure of
$\BA_{\bsymb{\a},\bsymb{\nu}}^{\bullet}$.
\end{thm}

\vskip0.2cm Our main goal here, as well as in the preceding paper
\cite{ESI}, is to study the spectrum of the operator
$\BA_{\bsymb{\a},\bsymb{\nu}}$ for different values of the
parameters $\a_\pm$. Informally, the mains results of both papers
can be summarized as follows: {\it the spectral properties of a
$K$-oscillator system can be described in terms of the corresponding
properties of $K$ appropriate one-oscillator systems}. To obtain
these one-oscillator systems, one divides the original graph into
$K$ pieces in such a way that each part contains only one point at
which an oscillator is attached, and these points should not
belong to the new boundary appearing as a result of the division.
On this new boundary we put an additional boundary condition, for
instance the Dirichlet condition. For our case ($\G=\R$ and the
oscillators attached at $\pm1$), it is most natural to take $x=0$
as the point of division. Let us denote the corresponding
operators by $\BA_{\R_\pm;\a_\pm;\nu_\pm}$; see \cite{ESI},
section 2.4 for details. \vskip0.2cm

The following theorem is proved in \cite{ESI}, Theorem 2.6.

\begin{thm}\label{ES1}
Let
\begin{equation*}
 \mu_{\pm} := \frac{\nu_{\pm}\sqrt2}{\a_{\pm}}.
\end{equation*}
\bigskip

\noindent { 1.} If $\mu_{\pm} > 1$, then $
\s_{a.c.}(\BA_{\bsymb{\a},\bsymb{\nu}})=[r_{0,0},\infty)=[(\nu_+^2+\nu_-^2)/2,
\infty).$

\bigskip

\noindent {2.} Let $\mu_+ =1 $ and $\mu_- >1$, or $\mu_- =1$ and
$\mu_+>1$. Then
\[ \s_{a.c.}(\BA_{\bsymb{\a},\bsymb{\nu}})=[\nu_-^2/2,\infty )\qquad {\rm{or}}\qquad
 \s_{a.c.}(\BA_{\bsymb{\a},\bsymb{\nu}})=[\nu_+^2/2,\infty)\]
respectively.

\bigskip

\noindent {3.} Let $\mu_+ = \mu_- =1$, then $
\s_{a.c.}(\BA_{\bsymb{\a},\bsymb{\nu}})=[0,\infty ).$

In all the cases {\rm 1 -- 3} the multiplicity function
$\gm_{a.c.}(\l;\BA_{\bsymb{\a},\bsymb{\nu}})$, is finite for all
$\l\in\s_{a.c.}(\BA_{\bsymb{\a},\bsymb{\nu}})$ and is given by
\begin{equation}\label{add3}\begin{split}
\gm_{a.c.}(\l;\BA_{\bsymb{\a},\bsymb{\nu}})=
&\sum_{n\in\N_0}\gm_{a.c.}(\l-\nu_-^2(n+1/2);\BA_{\R_+;\a_+;\nu_+})\\
+&\sum_{m\in\N_0}\gm_{a.c.}(\l-\nu_+^2(m+1/2);\BA_{\R_-;\a_-;\nu_-}).\end{split}
\end{equation}

\vskip0.2cm \noindent {4.} Let $max(\mu_+,\mu_-)<1$. Then
\[\s_{a.c.}(\BA_{\bsymb{\a},\bsymb{\nu}})=\R,\qquad
\gm_{a.c.}(\l;\BA_{\bsymb{\a},\bsymb{\nu}})\equiv\infty.\]
\end{thm}
\vskip0.2cm

In the present paper we are concerned with the point spectrum
below the threshold $r_{0,0}$ in the case that $\mu_+$ and $\mu_-$ are
both greater than 1. Below $N_-(\l;\BT)$, where $\l$ is a real
number, stands for the number of eigenvalues (counting
multiplicities) of a self-adjoint operator $\BT$, lying on the
half-line $(-\infty,\l)$, provided that this part of the spectrum
is discrete. We also set $N_+(\l;\BT)=N_-(-\l;-\BT)$.

\vskip0.2cm On the qualitative level, the main result of this
paper can be described as follows.

{\it For any $\mu_+,\mu_-<1$ the number
$N_-\left(r_{0,0};\BA_{\bsymb{\a},\bsymb{\nu}}\right)$ is finite and
asymptotically
\begin{equation}\label{2.inf}\begin{split}
N_-\left(r_{0,0};\BA_{\bsymb{\a},\bsymb{\nu}}\right) & \sim
N_-(\nu_+^2/2;\BA_{\R_+;\a_+;\nu_+})+N_-(\nu_-^2/2;\BA_{\R_-;\a_-;\nu_-}),\\
r_{0,0}=&(\nu_+^2+\nu_-^2)/2, \qquad \mu_\pm\downarrow 1.\end{split}
\end{equation}}

In order to give the precise formulation, we need to describe the
behaviour of the terms on the right-hand side of \eqref{2.inf},
and to explain what we mean when speaking about the asymptotics in
two parameters. To achieve the first goal, we present a result
which is a special case of Theorem 3.1 in \cite{S3}, see also
(3.10) in \cite{S2}. Let $\G=[a,b]$ (with the standard change if
$a=-\infty$ or $b=\infty$) be a finite or infinite interval and
$o\in\Int{\G}$. Consider the operator $\BA_{\G;\a;\nu}$ in
$L^2(\G\times\R)$, defined by the differential expression
\[\CA U=-U_{x^2}'' + \frac{\nu^2}{2} (-U_{q^2}''+  q^2 U)\]
and the matching condition
\[ U_x'(o+,q)-U_x'(o-,q)= \a qU(o,q),\]
cf \eqref{e1.1} and \eqref{e1.2}. If $\G\neq\R$, the Dirichlet or
the Neumann boundary condition is posed on $\p\G\times\R$. We do
not reflect the type of this condition in our notation. If
$\G=\R$, we drop the index $\G$ in the notation of the operator.
\begin{prop}\label{1.osc}
For any $\a\in(0,\nu\sqrt2)$ the spectrum of the operator
$\BA_{\a,\nu}$ below the point $\nu^2/2$ is non-empty and finite,
and the following asymptotic formula is satisfied:
\begin{equation}\label{0.6b}
N_-(\nu^2/2; \BA_{\G;\a;\nu}) \sim\frac{1}{4\sqrt{2(\mu-1)}},\qquad
\mu:=\frac{\nu\sqrt2}{\a}\downarrow 1.
\end{equation}
\end{prop}
It was assumed in \cite{S2} and \cite{S3}  that $\nu=1$, the
general case reduces to this special case by scaling. \vskip0.2cm

Our next theorem, together with the subsequent explanation of
uniformity of the asymptotics, gives the precise meaning to
\eqref{2.inf}. In the formulation of its second part an arbitrary
positive function $\psi(t)$ on $(0,1)$ which is $o(t^{-1/4})$ as $
t \to 0$ is involved. We also define the set
\begin{equation}\label{add4}
    \Omega_{\Psi} := \left \{ (x,y): \Psi(x) \le y \le 1,\ \ \ \Psi(y) \le x \le
    1 \right \},\ \ \ \Psi(t) = e^{-\psi(t)}.
\end{equation}
Note that the co-ordinate axes are tangents of infinite order to
$\Omega_{\Psi}$ at the origin.

\begin{thm} \label{ESII1}
\noindent {\rm {1.}} If $\mu_{\pm} := \sqrt 2 \nu_{\pm}/ \a_{\pm}
> 1$, then $\BA_{\bsymb{\a},\bsymb{\nu}}$ is bounded below and its
spectrum in $(-\infty, r_{0,0})$ is non-empty and finite.

\bigskip

\noindent{\rm{2.}} Let $\Psi$ be chosen as in \eqref{add4}.  Then,
uniformly for $(1-\mu_+^{-1},1-\mu_-^{-1}) \in \Omega_{\Psi}$,
\begin{equation}\label{1.10}
N_-\left(r_{0,0}; \BA_{\bsymb{\a},\bsymb{\nu}}\right) \sim
\frac{1}{4\sqrt{2(\mu_+ -1)}} \nonumber
     +   \frac{1}{4\sqrt{2(\mu_--1)}},\qquad \mu_\pm\downarrow 1.
\end{equation}
\end{thm}

Now, let us explain what we mean by `uniform asymptotics'. It means
that on the domain $(1-\mu_+^{-1},1-\mu_-^{-1})\in\Omega_{\Psi}$ there exists a
bounded function
$\Phi(\mu_+,\mu_-)$, such that $\Phi(\mu_+,\mu_-)\to 0$ as $\mu_\pm\to 1$ and
\begin{equation*}\begin{split}
    &\bigl|N_-\left(r_{0,0}; \BA_{\bsymb{\a},\bsymb{\nu}}\right)-
\frac{1}{4\sqrt2}((\mu_+ -1)^{-1/2}+ (\mu_--1)^{-1/2})\bigr|\\
&\le\Phi(\mu_+,\mu_-)\left((\mu_+ -1)^{-1/2}+ (\mu_--1)^{-1/2}\right).\end{split}
\end{equation*}

The technical ideas which lead to this result were explained in
the introduction to \cite{ESI}. Here we only note that for
$\mu_\pm\ge1$ the operator $\BA_{\bsymb{\a},\bsymb{\nu}}$ is
bounded below (see \thmref{qf.t}), which makes it possible to
apply the variational approach. In contrast to the operator domain
of the operator $\BA_{\bsymb{\a},\bsymb{\nu}}$, its quadratic form
domain for $\mu_\pm>1$ does not depend on the parameters $\a_\pm$.
This significantly simplifies the analysis. In particular, we do
not need to divide the graph into two parts, as we did in
\cite{ESI}; cf. (2.9) and Theorem 2.8 there. In our proof of
\thmref{ESII1} we will be dealing with the operators
$\BA_{\a_\pm;\nu_\pm}$ (i.e., the corresponding graph is $\G=\R$)
rather than with $\BA_{\R_\pm\a_\pm;\nu_\pm}$ as in \eqref{2.inf}.
According to Proposition \ref{1.osc}, the passage from $\R_\pm$ to
$\R$ does not affect the asymptotic behaviour of the function
$N_-$ for these operators.

\vskip0.2cm

We mostly use the same notation as in \cite{ESI}. However, in this
paper we have to take special care in order to distinguish between
the operators which correspond to the one-oscillator and to the
two-oscillator cases. We always denote the first as $\BA_{\a,\nu}$
and the second as $\BA_{\bsymb{\a},\bsymb{\nu}}$, with the
boldface $\a,\nu$ in the indices. Besides, we almost never drop
the index $\nu$ in the notation.

\section{Variational description of $\BA_{\bsymb{\a},\bsymb{\nu}}$ for
$\mu_{\pm} > 1$}
\label{bel} \subsection{The quadratic form
$\ba_{\bsymb{\a},\bsymb{\nu}}$.} If $U \sim \{u_{m,n}\} \in
\CD_{\bsymb\a,\bsymb\nu}$, the quadratic form
$\ba_{\bsymb{\a},\bsymb{\nu}}[U] :=
(\BA_{\bsymb{\a},\bsymb{\nu}}U,U)$ is given by
\begin{equation}\label{e1.6}
    \ba_{\bsymb{\a},\bsymb{\nu}}[U] = \ba[U] + \a_+ \bb_+[U] +
\a_- \bb_-[U],
\end{equation}
where, in the notation \eqref{sa.2},
\begin{gather}
  \ba[U]  = \sum_{m,n \in \mathbb{N}_0}
\int_\mathbb{R}\big( |u_{m,n}'(x)|^2
  + r_{m,n}|u_{m,n}|^2\big )dx, \label{1.7} \\
  \bb_+[U] =  \re\sum_{m,n \in \mathbb{N_0}}
  \sqrt {2m}\, u_{m,n}(1)
  \overline{u_{m-1,n}(1)}, \label{1.7x}\\
   \bb_-[U] =  \re\sum_{m,n \in {\N_0}}\sqrt {2n}\, u_{m,n}(-1)
   \overline{u_{m,n-1}(-1)}\label{1.7y}.
\end{gather}
In \eqref{1.7x} and  \eqref{1.7y} we took by default that
$u_{-1,n}\equiv 0$ and $u_{m,-1}\equiv 0$ for all $m,n\in\N$.

The quadratic form $\ba$ (which is the same as
$\ba_{\bsymb{0},\bsymb{\nu}}$) is positive definite in $\GH$.
Completing the set $\CD_{\bsymb{0},\bsymb{\nu}}$ with respect to
the `energy metric' $\ba[U]$, we obtain a Hilbert space which we
denote by $\bd$.

\vskip0.2cm Let us define $H^1_\g$, where $\g>0$ is a real
parameter, to be the Sobolev space $H^1(\R)$ with the scalar
product
\begin{equation}\label{sob}
(u_1,u_2)_\g=\int_\R\left(u'_1(x)\overline{u'_2(x)}+\g^2u_1(x)
\overline{u_2(x)}\right)dx
\end{equation}
and the corresponding norm $\|u\|_\g$. The space $\bd$ can be
naturally identified with the orthogonal sum of the spaces
$H^1_{\sqrt{r_{m,n}}}$. The topology in $\bd$ does not depend on
the values of $\nu_\pm.$

Our next goal is to prove the following
\begin{thm}\label{qf.t}
Let $\mu_\pm\ge1$. Then the quadratic form
$\ba_{\bsymb{\a},\bsymb{\nu}}$ is bounded below. If $\mu_\pm>1$,
then $\ba_{\bsymb{\a},\bsymb{\nu}}$ is closed on $\bd$ and the
corresponding self-adjoint operator in $\GH$ coincides with
$\BA_{\bsymb{\a},\bsymb{\nu}}$.
\end{thm}

For the proof we need some auxiliary material. Let $\CF_\g$  be
the two-dimensional space of functions $v\in H^1_\g$ which for
$x\neq\pm1$ satisfy the equation \[-v''+\g^2 v=0.\] Evidently,
each function $v\in H^1_\g$ is uniquely determined by its values
at the points $\pm1$. The space $\CF_\g$ was discussed in
\cite{ESI}, sec. 3.1. In particular, it was shown there that for
any $v\in \CF_\g$ one has
\begin{equation}\label{1c}
[v'](p)=-\frac{2\g}{1-e^{-4\g}}\left(v(p)-e^{-2\g}v(-p)\right),
\qquad p=\pm1.
\end{equation}
It follows from \eqref{1c} that the mapping $v\mapsto
([v'](1),[v'](-1))$ maps $\CF_\g$ onto $\C^2$.

Denote by $\Pi_\g$ the operator of ortogonal projection (in the scalar
product \eqref{sob}) of the space $H_\g^1$ onto $\CF_\g$.
\begin{lem}\label{l1}
For any $u\in H^1_\g$ its projection $\Pi_\g u$ is the function
$v\in\CF_\g$, defined by the conditions
\begin{equation}\label{1d}
v(\pm1)=u(\pm1).
\end{equation}
\end{lem}
\begin{proof}
Let $v,w\in\CF_\g$. We have
\[\begin{split}
&(u-v,w)_\g=\left(\int_{-\infty}^{-1}+\int_{-1}^1+\int_1^\infty\right)
\left((u'-v')\overline{w'}+\g^2(u-v)\overline{w}\right)dx\\
&=\left(\int_{-\infty}^{-1}+\int_{-1}^1+\int_1^\infty\right)
(u-v)(\overline{-w''+\g^2{w}})dx\\
&-(u(1)-v(1))\overline{[w'](1)}
-(u(-1)-v(-1))\overline{[w'](-1)}.\end{split}\] The integrand in
the second line vanishes and we get
\[(u-v,w)_\g=-(u(1)-v(1))\overline{[w'](1)}-
(u(-1)-v(-1))\overline{[w'](-1)}.\] By \eqref{1c}, the set of all
possible pairs $([w'](1),[w'](-1))$ covers the whole of $\C^2$
which implies the result.
\end{proof}

\vskip0.2cm
\begin{lem}\label{l2}
For all $u\in H^1_{\g},$
\begin{equation}\label{2.1}
    2\g(|u(-1)|^2+|u(1)|^2) \le (1+e^{-2\g}) \int_
    \R \left(|u'|^2+\g^2|u|^2\right)dx.
\end{equation}
The constant is optimal. The equality in \eqref{2.1} is attained
on the one-dimensional subspace in $H^1_\g$ formed by the
functions $v\in\CF_\g$ such that $v(1)=v(-1)$.
\end{lem}
\begin{proof} Given a function $u\in H^1_\g$, take $v=\Pi_\g u$.
Then
\[ \|u-v\|^2_\g=\|u\|^2_\g-(v,u)_\g=\|u\|^2_\g-
\int_\R\left(v'\overline{u'}+\g^2v\overline{u}\right)dx.\]
Integrating by parts as in Lemma \ref{l1} and denoting $u(1)=A,\
u(-1)=B$, we get
\begin{gather*}
\|u-v\|_\g^2= \|u\|_\g^2+\overline{A}[v'](1)+\overline{B}[v'](-1).
\end{gather*}
On using \eqref{1c} and \eqref{1d}, we find from here:
\begin{gather*}
0\le\|u-v\|_\g^2=\|u\|_\g^2-\frac{2\g}{1-e^{-4\g}}\left(\overline{A}
(A-e^{-2\g}B)+\overline{B}(B-e^{-2\g}A)\right)\\
=\|u\|_\g^2-\frac{2\g}{1+e^{-2\g}}(|A|^2+|B|^2)-\frac{2\g
e^{-2\g}} {1-e^{-4\g}}|A-B|^2,
\end{gather*}
whence the Lemma.
\end{proof}

\subsection{Proof of \thmref{qf.t}.}\label{qft}
We obtain from \eqref{1.7x}:
\begin{equation*}
 \bb_+[U] \le \frac{1}{2} \sum_{m,n \in \mathbb{N}_0} \left(\sqrt{2 m} + \sqrt{2(m+1)}
  \right)|u_{m,n}(1)|^2
   \le   \sum_{m \in \mathbb{N}_0,n\in\N} \sqrt{2m + 1}|u_{m,n}(1)|^2
  \end{equation*}
  and similarly
\begin{equation*}
  \bb_-[U] \le
      \sum_{m\in\N,,n \in \mathbb{N}_0} \sqrt{2n +1}|u_{m,n}(-1)|^2.
  \end{equation*}

Given a number $k\ge -r_{0,0}$, denote
\begin{equation}\label{1.8x}
 \g_{m,n}(k) = \sqrt{r_{m,n} +k}.
\end{equation}
The conditions $\mu_+,\mu_-\ge1$ imply
\[
\max(\a_+\sqrt{2m + 1},\a_-\sqrt{2n + 1})\le 2\g_{m,n}(0).\]
 Hence,
\[\begin{split}
&\a_+\sqrt{2m + 1}|u_{m,n}(1)|^2+\a_-\sqrt{2n +1}|u_{m,n}(-1)|^2\\
&\le 2\g_{m,n}(0)\left(|u_{m,n}(1)|^2+|u_{m,n}(-1)|^2\right).
\end{split}\]
Applying Lemma \ref{l2} with $\g = \g_{m,n}(k)$ and $k$ a positive
constant to be chosen later, we obtain
\begin{equation}\label{1.8a}
\a_+\sqrt{2m + 1}|u_{m,n}(1)|^2+\a_-\sqrt{2n +1}|u_{m,n}(-1)|^2\le
C(m,n,k)\|u_{m,n}\|^2_{H^1_{\g_{m,n}(k)}}
\end{equation}
where
\[C(m,n,k)=
\frac{\g_{m,n}(0)}{\g_{m,n}(k)}(1+e^{-2\g_{m,n}(k)}).\]

Now we show that
\begin{equation}\label{1.8b}
C(m,n,k)\le 1,\qquad \forall m,n\in\N_0,
\end{equation}
provided that $k$ is large enough. To this end, consider the function
\[f_k(t)=(1-kt^{-2})^{1/2}\left(1+e^{-2t}\right),\qquad t\ge k^{1/2},\ k>0,\]
then
\[C(m,n,k)=f_k(\g_{m,n}(k)).\]
Note that $f_k(k^{1/2})=0$ and $f_k(t)\to 1$ as $t\to\infty$. Hence,
\eqref{1.8b}  will be proven if we show that
$f'_k(t)\ge 0$ for all $t$.

We have
\[\begin{split} &f'(t)=\frac{k(1+e^{-2t})}{t^3(1-kt^{-2})^{1/2}}
-2(1-kt^{-2})^{1/2}e^{-2t}\\
&=\frac{(1+e^{-2t}+2te^{-2t})k-2t^3e^{-2t}}{t^3(1-kt^{-2})^{1/2}}
\ge \frac{k-2t^3e^{-2t}}{t^3(1-kt^{-2})^{1/2}},\end{split}\]
 and the desired result follows for $k\ge 27
e^{-3}/4=\max\left(2t^3e^{-2t} \right)$.

On taking $k$ such that \eqref{1.8b} is satisfied, we derive from
\eqref{1.8a}:
\[ |\a_+\bb_+[U]+\a_-\bb_-[U]|\le \sum_{m,n\in\N_0}\|u_{m,n}\|^2_{H^1_{\g_{m,n}(k)}}
=\ba[U]+k\|U\|^2_\GH.\]
 So, the boundedness below of $\ba_{\bsymb{\a},\bsymb{\nu}}$ for
 all $\mu_\pm\ge1$ is established. The closedness of this
 quadratic form for all $\mu_\pm>1$ easily follows from here,
 cf. \cite{BS}. Since the operator
$\BA_{\bsymb{\a},\bsymb{\nu}}$ has a unique self-adjoint
realization, it necessarily coincides with the operator associated
with the quadratic form $\ba_{\bsymb{\a},\bsymb{\nu}}$.

\section{The spectrum of $\BA_{\bsymb{\a},\bsymb{\nu}}$ below
$r_{0,0}$.}\label{spb}

We next prove that the spectrum below $r_{0,0}$ is finite and
non-empty, and in the process, give an alternative proof of part 1
of Theorem \ref{ESII1}. Our argument is similar to the one in \cite{S3}
where the one-oscillator case was studied.
\subsection{Finiteness.}
For some $L \in \mathbb{N},$ let us consider the quadratic form
$\ba_{\bsymb{\a},\bsymb{\nu}}$, see \eqref{e1.6},  on the set
\begin{equation}\label{3.1}
\bd^{(L)}=\bigl\{U\sim\{u_{m,n}\}: u_{m,n}(\pm1)=0,\ m+n\le
L\bigr\}.
\end{equation}
For the operator $\BA_{\bsymb{\a},\bsymb{\nu}}^{(L)}$, associated
with the quadratic form
$\ba_{\bsymb{\a},\bsymb{\nu}}\res\bd^{(L)}$, the subspace
\[\GH^{(L)}=\bigl\{U\sim\{u_{m,n}\}: u_{m,n}\equiv 0,\ m+n>
L\bigr\}\] is invariant, and the part
$\BA_{\bsymb{\a},\bsymb{\nu}}^{(L,-)}$ of
$\BA_{\bsymb{\a},\bsymb{\nu}}^{(L)}$ in $\GH^{(L)}$ decomposes in
the orthogonal sum:
\begin{equation}\label{3.9}
\BA_{\bsymb{\a},\bsymb{\nu}}^{(L,-)}= {\sum_{m+n \le L}}^{\oplus}
\left(\BA+r_{m,n}\right),
\end{equation}
where $\BA=-\frac{d^2}{dx^2} $ with domain $H^2( \mathbb{R})$.
Since
\[
\sigma(\BA) = \sigma_{a.c.}(\BA) = [0,\infty),
\]
it follows that
\begin{equation}\label{3.10}
    \sigma(\BA_{\bsymb{\a},\bsymb{\nu}}^{(L,-)}) =
    \sigma_{a.c.}(\BA_{\bsymb{\a},\bsymb{\nu}}^{(L,-)})
    =[r_{0,0},\infty).
\end{equation}
An explicit expression for the multiplicity function
$\gm_{a.c.}(\l;\BA_{\bsymb{\a},\bsymb{\nu}}^{(L,-)})$ immediately
follows from \eqref{3.9}, but this is omitted.

On repeating the argument in section \ref{qft} with $k = 0$, we
have that for $U\in\bd,\ U\perp\GH^{(L)}$
\begin{eqnarray*}
   \ba_{\bsymb{\a},\bsymb{\nu}}[U]&\ge & \sum_{m+n>L}
   \left(1- \max\{\mu_+^{-1},\mu_-^{-1}\}
   (1+e^{-2\gamma_{m,n}(0)})\right) \\
   &\times & \int_{ \mathbb{R}}\left(|u_{m,n}'|^2
   +r_{m,n}|u_{m,n}|^2\right)dx .
\end{eqnarray*}

Let $\BA_{\bsymb{\a},\bsymb{\nu}}^{(L,+)}$ stand for the part of
$\BA_{\bsymb{\a},\bsymb{\nu}}$ in the subspace
$(\GH^{(L)})^\perp$. It follows from the above inequality that for
any $\l_0>0$ it is possible to choose $L$ sufficiently large, to
ensure that
\begin{equation}\label{3.11}
    (\BA_{\bsymb{\a},\bsymb{\nu}}^{(L,+)}U,U)\ge \l_0\|U\|^2.
\end{equation}
Hence, in view of (\ref{3.10}),
\begin{equation}\label{3.12}
  \sigma(\BA_{\bsymb{\a},\bsymb{\nu}}^{(L)})= [r_{0,0},\infty),
\end{equation}
\begin{equation}\label{3.12a}
  \sigma_{a.c.}(\BA_{\bsymb{\a},\bsymb{\nu}}^{(L)} ) \supseteq
  [r_{0,0},\lambda_0).
\end{equation}

The passage from the operator $\BA_{\bsymb{\a},\bsymb{\nu}}$ to
$\BA_{\bsymb{\a},\bsymb{\nu}}^{(L)}$ corresponds to the passage
from the quadratic form domain $\bd$ to its subspace $\bd^{(L)}$ of
finite co-dimension. In its turn, this corresponds to a finite
rank perturbation of the resolvent. Such perturbations do not
affect the absolutely continuous spectrum and its multiplicity.
Hence,
\[\gm_{a.c.}(\l;\BA_{\bsymb{\a},\bsymb{\nu}}
)=\gm_{a.c.}(\l;\BA^{(L)}_{\bsymb{\a},\bsymb{\nu}} ),\qquad
\l\in[r_{0,0},\infty).\] This immediately leads to
\eqref{add3} for $r_{0,0}\le\l<\l_0$ and therefore, for all $\l\ge
r_{0,0}$.

Besides, the number of eigenvalues of
$\BA_{\bsymb{\a},\bsymb{\nu}}$ which may appear below $r_{0,0}$
under such a  perturbation, does not exceed the rank of the perturbation
and hence, is finite.

\subsection{Non-emptiness of $\s_p$.}
To prove that the spectrum below $r_{0,0}$ is non-empty, and hence
complete the proof of \thmref{ESII1}, we apply the argument used
to prove the analogous result in \cite{S1}, Theorem 6.2. It is
sufficient to find a function $U \in \bd$ which is such that
\begin{equation}\label{1.9}
    \ba_{\bsymb{\a},\bsymb{\nu}}[U] < r_{0,0}\|U\|^2_{\GH}.
\end{equation}
Choose $U\sim\{u_{m,n}\}$ as follows. We take
\[ u_{0,0}(x)=-\vare^{-1/2}\min(1,e^{-(\vare|x|-1)}),\]
with $\vare\in(0,1)$ to be chosen later. Note that
\[\int_\R|u'_{0,0}|^2dx=1,\qquad u_{0,0}(\pm1)=-\vare^{-1/2}.\]
We also take $u_{1,0}(x)=e^{-|x-1|},\ u_{0,1}(x)=e^{-|x+1|}$, then
$u_{1,0}(1)=u_{0,1}(-1)=1$ and
\[\int_\R|u'_{1,0}|^2dx=\int_\R|u'_{0,1}|^2dx=\int_\R|u_{1,0}|^2dx=\int_\R|u_{0,1}|^2dx=1.\]
We take all the other components $u_{m,n}$ to be zero.
 For such $U$ we have
\begin{eqnarray*}
    &\ba_{\bsymb{\a},\bsymb{\nu}}[U] -
    r_{0,0}\|U\|^2_{\GH}\\
    &= \int_{\R}\left(|u'_{0,0}|^2 + |u'_{1,0}|^2 +
    |u'_{0,1}|^2+\nu_+^2|u_{1,0}|^2+
    \nu_-^2|u_{0,1}|^2\right)dx \\
    &+ \sqrt 2
    \a_+ u_{1,0}(1)u_{0,0}(1) + \sqrt 2
    \a_-u_{0,1}(-1) u_{0,0}(-1)\\
    &=3+   \nu_+^2+
    \nu_-^2 - \vare^{-1/2}\sqrt 2(\a_+ +\a_-).
\end{eqnarray*}
On choosing $\vare $ sufficiently small we obtain a function $U$
which satisfies (\ref{1.9}).
 This
completes the proof of part 1 of Theorem 2.1.

\section{Asymptotics: reduction to a problem in
$\ell^2$}\label{asym} \subsection{Removing the component
$u_{0,0}$.} In what follows it is convenient for us to consider
the quadratic form $\ba_{\bsymb{\a},\bsymb{\nu}}$, defined in \eqref{e1.6}, for the
elements $U\sim\{u_{m,n}\}\in\bd$ subject to the additional
conditions
\begin{equation}\label{x1}
u_{0,0}(1)=u_{0,0}(-1)=0.
\end{equation}
 For any $\a_\pm<\nu_\pm\sqrt2$ the quadratic form
$\ba_{\bsymb{\a},\bsymb{\nu}}$, restricted to this domain,
generates in $\GH$ a self-adjoint operator, for which the subspace
\[\GH_{0,0}=\left\{U\sim \{u_{0,0},0,0,\ldots\}\right\}\]
is invariant. The part of this operator in $\GH_{0,0}$ is
$-u''_{0,0}+r_{0,0}u_{0,0}$ under the conditions \eqref{x1} and it
has no spectrum below $r_{0,0}$. Removing this subspace yields
the Hilbert space
\[\GH^\circ=\left\{U\sim\{u_{m,n}\}:u_{0,0}\equiv0\right\}\]
and the quadratic form
$\ba^\circ_{\bsymb{\a},\bsymb{\nu}}=\ba_{\bsymb{\a},\bsymb{\nu}}\res\bd^\circ$.

Below we denote
\[\g_{m,n}=\g_{m,n}(-r_{0,0})=\sqrt{\nu_+^2m+\nu_-^2n},\]
cf \eqref{1.8x}. We shall consider $\bd^\circ$ as a Hilbert space
with the norm given by
\begin{equation}\label{x2}
\|U\|^2_{\bd^\circ}=\ba^\circ[U]-r_{0,0}\|U\|^2_{\GH^\circ}
=\sum_{m+n>0}\int_\R\left(|u_{m,n}'|^2+\g_{m,n}^2|u_{m,n}|^2\right)dx
\end{equation}
and the corresponding scalar product $(.,.)_{\bd^\circ}$. The norm
$\|U\|_{\bd^\circ}$ and the ``energy norm'' $\sqrt{\ba[U]}$ are
equivalent on $\bd^\circ$. On the whole of $\bd$ this is not true.
This explains, why the passage from $\bd$ to
$\bd^\circ$ is useful.

Let $\BA^\circ_{\bsymb{\a},\bsymb{\nu}}$ stand for the
self-adjoint operator in $\GH^\circ$, associated with the
quadratic form $\ba^\circ_{\bsymb{\a},\bsymb{\nu}}$. It follows
from the variational argument that
\begin{equation}\label{x3}
0\le
N_-(r_{0,0};\BA_{\bsymb{\a},\bsymb{\nu}})-N_-(r_{0,0};\BA^\circ_{\bsymb{\a},\bsymb{\nu}})\le
2,\qquad \forall \mu_\pm>1.
\end{equation}
Therefore, both counting functions have the same asymptotic
behaviour as $\mu_\pm\downarrow 1$.

\vskip0.2cm

According to the variational principle,
\begin{equation}\label{x4}
N_-(r_{0,0};\BA^\circ_{\bsymb{\a},\bsymb{\nu}})=\min_{E\in\CE}\codim
E
\end{equation}
where $\CE$ is the set of all subspaces $E\subset\bd^\circ$ such
that
\[ \ba^\circ_{\bsymb{\a},\bsymb{\nu}}[U]\ge r_{0,0}\|U\|^2_{\GH^\circ},
\qquad\forall U\in E.\] The latter inequality can be re-written as
\begin{equation}\label{x5}
\|U\|^2_{\bd^\circ}+\a_+\bb_+[U]+\a_-\bb_-[U]\ge 0,\qquad\forall
U\in E.
\end{equation}
\subsection{Shrinking the space.}
Our next goal is to show that it is enough to take the maximum in
\eqref{x4} over the set of subspaces $E\subset\CF$ where
\begin{equation*}
\CF={\sum_{m+n>0}}^\oplus\CF_{\g_{m,n}}
\end{equation*}
(recall that the two-dimensional spaces $\CF_\g$ were defined in section 2.1).
Indeed, in the variational description of the non-zero spectrum of
a self-adjoint operator $\BT$ one can always consider only the
subspaces orthogonal to $\ker \BT$. Let us apply this remark to
the operator $\BB$ in $\bd^\circ$, generated by the right-hand
side in \eqref{x5}. It follows from Lemma \ref{l1} that the
orthogonal complement to $\CF$ in $\bd^\circ$ is given by
\[\CF^\perp=\left\{U\sim\{u_{m,n}\}:u_{m,n}(1)= u_{m,n}(-1)=0.\right\}\]
Therefore, $\CF^\perp\subset\ker\BB$, which yields the desired
result; see \cite{S2}, proof of Theorem 3.1, or \cite{S3}, proof
of Theorem 10.1, for further details.

\vskip0.2cm

Now we construct a convenient orthogonal basis in $\CF$. It is enough
to choose a basis in each component $\CF_{\g_{m,n}}$. To simplify
notation, in calculations below we drop the indices $m,n$.

The functions $u^\pm(x)=e^{-\g|x\mp1|}$ form a linear basis in
$\CF_{\g}$. We have
\[\|u^\pm\|^2_\g=2\g, \qquad (u^+,u^-)_\g=2\g e^{-2\g}\]
where the norm and the scalar product are taken in $H^1_\g$, see
\eqref{sob}. Let now
\begin{equation}\label{4.7}
    v^+:= \frac{u^++\varkappa u^-}{\|u^+ +\varkappa
    u^-\|_{\gamma}},\ \ \ v^-:= \frac{u^- +\varkappa u^+}{\|u^- +\varkappa
    u^+\|_{\gamma}},
\end{equation}
for a constant $\varkappa.$  These are normalized and are
orthogonal in $H^1_{\gamma}$ if and only if
$\varkappa^2+2e^{2\gamma}\varkappa +1 =0.$ We choose the root
\begin{equation}\label{4.7a}
\vark=-e^{2\gamma}+\sqrt{e^{4\gamma}-1}=-\frac1{2}e^{-2\gamma}(1+O(e^{-4\gamma}));
\end{equation}
then
\begin{equation}\label{4.7b}
\rho^2:=\|u^\pm +\varkappa u^\mp\|^2_{\gamma}=
2\g(1+\vark^2+2\vark e^{-2\gamma})=2\g(1+O(e^{-4\gamma})).
\end{equation}
 Also, using the equation for $\vark$, we find
\[ v^+(1)=v^-(-1)=\wh{\rho}\,^{-1},\qquad v^+(-1)=v^-(1)=-\vark\wh{\rho}\,^{-1}\]
where
\[\wh{\rho}=\rho(1-e^{-4\g})^{-1/2}.\]

\vskip0.2cm Below we indicate the dependence of $\g,\ \varkappa$
and $\rho$ on $m,n$. In particular, we write $v^\pm_{m,n}$. Note
that by \eqref{4.7a}, \eqref{4.7b} we have
\begin{equation}\label{4.7c}
\varkappa_{m,n} =
  -\frac{e^{-2\gamma_{m,n}}}{2}\left(1+O(e^{-4\gamma_{m,n}})\right),\qquad
    \rho_{m,n} = \sqrt{2\gamma_{m,n}}
  \left(1+O(e^{-4 \gamma_{m,n}})\right) .\end{equation}

Let $U\sim\{C^+_{m,n}v^+_{m,n}+ C^-_{m,n}v^-_{m,n}\}\in\CF$, then
the mapping \[U\mapsto\CC=\{C^+_{m,n},C^-_{m,n}\}\] is an isometry
of $\CF$ onto the Hilbert space
$\CG=\ell^2(\N_0^2\setminus\{(0,0)\})$. We denote by $\CG^\pm$ the
subspaces in $\CG$, formed by the elements
\[ \CC^+=\{C^+_{m,n},0\},\qquad \CC^-=\{0,C^-_{m,n}\}\]
respectively. On $\bd^\circ$ the quadratic forms $\bb_\pm$ become
\[\begin{split}
  &\bb_+[U] = \bb'_+[\CC]\\ =\sum_{m+n>0} \frac{\sqrt{2m}}{\wh{\rho}_{m,n}\wh{\rho}_{m-1,n}}
  &\re[(C^+_{m,n}- \varkappa_{m,n}C^-_{m,n})
   \overline{(C^+_{m-1,n}- \varkappa_{m-1,n}C^-_{m-1,n})}],\\
  &\bb_-[U] =\bb'_-[\CC]\\=\sum_{m+n>0} \frac{\sqrt{2n}}{\wh{\rho}_{m,n}\wh{\rho}_{m,n-1}}
  &\re[(C^-_{m,n}- \varkappa_{m,n}C^+_{m,n})
   \overline{(C^-_{m,n-1}-\varkappa_{m,n-1}C^+_{m,n-1})}],\end{split}\]
and the quadratic form $\ba^\circ_{\bsymb{\a},\bsymb{\nu}}$
becomes
\begin{equation*}
    \ba'_{\bsymb{\a},\bsymb{\nu}}[\CC]=\|\CC\|^2_{\CG}+\a_+\bb'_+[\CC]+\a_-\bb'_-[\CC].
\end{equation*}
Denote by $\BB'_\pm$ the operators in $\CG$ associated with the
quadratic forms $\bb'_\pm$; then the operator associated with
$\ba'_{\bsymb{\a},\bsymb{\nu}}$ is $\BI+\a_+\BB'_+ +\a_-\BB'_-$.

It follows from this construction and \eqref{x4}, \eqref{x5} that
\begin{equation}\label{4.10}
    N_-(r_{0,0};\BA^\circ_{\bsymb{\a};\bsymb{\nu}})=N_-(0;\BI+\a_+\BB'_+
    +\a_-\BB'_-)=N_+(1;-\a_+\BB'_+ -\a_-\BB'_-).
\end{equation}

Consider now the case when one of the parameters $\a_\pm$ is equal
to zero. Below we denote
\[\bsymb{\a}_+=(\a_+,0),\qquad \bsymb{\a}_-=(0,\a_-).\]
For $\bsymb{\a}=\bsymb{\a}_\pm$ the equality \eqref{4.10} can be
re-written in the standard form of the Birman -- Schwinger
principle:
\begin{equation}\label{4.11}
    N_-(r_{0,0};\BA^\circ_{\bsymb{\a}_+;\bsymb{\nu}})=N_+(\a_+^{-1};-\BB'_+),\qquad
     N_-(r_{0,0};\BA^\circ_{\bsymb{\a}_-;\bsymb{\nu}})=N_+(\a_-^{-1};-\BB'_-).
\end{equation}

\subsection{Structure of the operators
$\BB'_\pm$.} Denote by $\bb''_\pm$ the leading terms in the
expressions for $\bb'_\pm$, i.e.
\[\begin{split}
&\bb_+''[\CC] = \bb_+''[\CC^+]=\sum_{m+n>0}
\frac{\sqrt{2m}}{\wh{\rho}_{m,n}\wh{\rho}_{m-1,n}}
  {\re}(C^+_{m,n} \overline{C^+_{m-1,n}}), \\
&  \bb_-''[\CC] =\bb_-''[\CC^-]=\sum_{m+n>0}
\frac{\sqrt{2n}}{\wh{\rho}_{m,n}\wh{\rho}_{m,n-1}}
  \re(C^-_{m,n}
   \overline{C^-_{m,n-1}}).\end{split}\]
Let $\BB_\pm''$ stand for the corresponding operators in
$\CG^\pm$.

Now we are in a position to explain the scheme of our further
analysis. It is natural to expect that the number
$N_-(0;\BI+\a_+\BB'_+ +\a_-\BB'_-)$, is close to
$N_-(0;(\BI_+ +\a_+\BB''_+)\oplus (\BI_- +\a_-\BB''_-))$. Indeed,
consider the operator
\begin{equation}\label{4.7x}\begin{split}
&\BX_{\bsymb{\a}}:=(\BI+\a_+\BB'_+ +\a_-\BB'_-)-
(\BI_++\a_+\BB''_+)\oplus (\BI_-+\a_- \BB''_-)\\ &=
\a_+(\BB'_+ -(\BB''_+\oplus 0))+\a_-(\BB'_- -(0\oplus\BB''_-)),
\end{split}\end{equation}
then
\[(\BX_{\bsymb{\a}}\CC,\CC)_\CG=\a_+(\bb'_+[\CC]-\bb''_+[\CC^+])
+\a_-(\bb'_-[\CC]-\bb''_-[\CC^-]).\]
This
quadratic form is expressed by a sum of terms
with exponentially decaying coefficients, and adding this sum
cannot affect the asymptotic behaviour of the function $N_-$.
Further, the behaviour of $N_-$ for the operator involving $\BB''_\pm$ is easy to
understand, due to its special structure.

So, our immediate task is to take care of the errors coming
from the difference $\bb'_\pm[\CC]-\bb''_\pm[\CC]$. Each term in
these quadratic forms involves at least one of the factors
$\vark_{m,n}, \vark_{m-1,n}, \vark_{m,n-1}$.  We have
\begin{equation*}
 \g_{m,n}=\sqrt{\nu_+^2m +\nu_-^2 n}\ge\d'\sqrt{m+n},\qquad \d'=\min(\nu_+,\nu_-).
\end{equation*}
Note also that by \eqref{4.7c} the factors
$\sqrt{2m}(\wh{\rho}_{m,n}\wh{\rho}_{m-1,n})^{-1}$ and
$\sqrt{2n}(\wh{\rho}_{m,n}\wh{\rho}_{m,n-1})^{-1}$ appearing in the
expressions for $\bb'_\pm,\bb''_\pm$ are bounded uniformly in
$m,n$. Taking this into account, applying the Cauchy -- Schwartz
inequality, and using the asymptotic result \eqref{4.7c} for
$\vark_{m,n}$, we come to the inequality
\[|\bb'_\pm[\CC]-\bb''_\pm[\CC]|\le
c\sum_{m+n>0}e^{-2\d\sqrt{m+n}}|C_{m,n}^\pm|^2,\] with some
$c<\infty$ and a positive $\d<\d'$. Now it follows from the
variational principle that the consecutive eigenvalues of the
operator $|\BX_{\bsymb{\a}}|$
do not exceed the numbers $c\max(\a_+,\a_-)e^{-2\d\sqrt{m+n}}$,
repeated twice and then rearranged in decreasing order.
Hence, given an $\vare>0$, we derive an estimate, uniform in
$\a_\pm\le \nu_\pm\sqrt2$:
\begin{equation}\label{4.21}
    N_+(\vare;|\BX_{\bsymb{\a}}|)\le \#\left\{(m,n)\in\N^2:C_0
e^{-2\d\sqrt{m+n}}>\vare\right\}\le R\log^4(K/\vare),
\end{equation}
with some $R,K>0$. Note that another way to obtain this inequality
is based on the connection between the eigenvalues and the
approximation numbers of a compact operator, see \cite{EE}.

\vskip0.2cm

Now, let us consider the operator $\BA_{\bsymb{\a}_+;\bsymb\nu}$.
Since $\a_-=0$, the variable $q_-$ can be separated and the
operator decomposes into the orthogonal sum (see (1.5) in
\cite{ESI})
\begin{equation*}
\BA_{\bsymb{\a}_+;\bsymb\nu}
={\sum_{n\in\N_0}}^\oplus(\BA_{\a_+;\nu_+}+\nu_-^2(n+1/2)).
\end{equation*}
This decomposition yields
\begin{equation}\label{1.6a}
N_-(r_{0,0}; \BA_{\bsymb{\a}_+;\bsymb\nu})
=\sum_{n\in\N_0}N_-(\nu_+^2/2-\nu_-^2n; \BA_{\a_+;\nu_+}).
\end{equation}
For $\a_+\le\nu_+\sqrt2$ the operator $\BA_{\a_+;\nu_+}$ is
non-negative (see \cite{S1}), therefore the sum in \eqref{1.6a}
has only a finite number of non-zero terms. Besides, the terms
corresponding to any $n>0$, are finite, since the essential
spectrum of $\BA_{\a_+,\nu_+}$ is $[\nu_+^2,\infty)$. Taking into
account that by \eqref{x3} the asymptotic behaviour as
$\a_+\to\nu_+\sqrt2$ of the function $N_-(r_{0,0}; .)$ for the
operators $\BA_{\bsymb{\a}_+;\bsymb\nu}$ and
$\BA^\circ_{\bsymb{\a}_+;\bsymb\nu}$ is the same, we conclude from
\eqref{0.6b} that \[N_-(r_{0,0};
\BA^\circ_{\bsymb{\a}_+;\bsymb\nu})\sim N_-(\nu_+^2/2;
\BA_{\a_+;\nu_+}) \sim\frac{1}{4\sqrt{2(\mu_+ -1)}},\qquad
\mu_+\downarrow 1.\] From the last equality and \eqref{4.11} we
derive that
\begin{equation*}
 N_-(0;\BI+\a_+\BB'_+)= N_+(\a_\pm^{-1};-\BB'_+)\sim\frac{1}{4\sqrt{2(\mu_+
-1)}},\qquad \mu_\pm\downarrow1.
\end{equation*}
The analogous equality is valid for the operator $\BB'_-$.

The same asymptotic formula holds for the operators $\BB''_\pm$:
\begin{equation}\label{1.6b}
 N_-(0;\BI_\pm+\a_\pm\BB''_\pm)\sim\frac{1}{4\sqrt{2(\mu_+
-1)}},\qquad \mu_\pm\downarrow1.
\end{equation}
This follows (for the `plus' sign, say) from the evident equality
\[N_-(0;\BI_+ +\a_+\BB''_+)=N_-(0;\BI +\a_+\BB''_+\oplus 0)\]
and from the estimate \eqref{4.21} for the case $\a_-=0$.

\section{ Proof of \thmref{ESII1}, part 2}\label{finish}
The proof is based upon \eqref{4.10} and the equality
\begin{equation*}
    \BI+\a_+\BB'_+ +\a_-\BB'_- =(\BI_+ +\a_+\BB''_+)\oplus(\BI_-
    +\a_-\BB''_-)+\BX_{\bsymb{\a}}
\end{equation*}
where the last term is given by \eqref{4.7x}.

Set $ \y_{\pm} := \mu_{\pm} - 1 =\frac{\nu_{\pm}\sqrt2}{\a_\pm}-1$
and $M = (4\sqrt2)^{-1}$. Then \eqref{1.6b} means that there exist
two non-negative functions $\varf_\pm(\mu_\pm)$, defined for
$\mu_\pm>1$, vanishing as $\mu_\pm\to 1$ and such that
\begin{equation}\label{5.2}
\left|N_-(0;\BI_\pm + \a_{\pm} \BB''_{\pm} )-M
(\mu_\pm-1)^{-1/2}\right| \le\varf_\pm (\mu_\pm)(\mu_\pm-1)^{-1/2} .
\end{equation}
To determine the asymptotic behaviour of $N_-(0;\BI+\a_+\BB'_+
+\a_-\BB'_-)$, we have to estimate the smallest co-dimension of
subspaces in $\CG$ on which
\begin{equation}\label{5.2a}
     \|\CC^+\|_\CG^2+\|\CC^-\|_\CG^2+\a_+\bb''[\CC^+]+
\a_-\bb''[\CC^-]+(\BX_{\bsymb{\a}}\CC,\CC)_\CG\ge 0
\end{equation}
for all $\CC$. By (\ref{4.21}), for any $\vare>0$ there
exists a subspace $ \CK(\vare) \subset \CG $ such that
\begin{equation}\label{4}
\codim \CK(\vare)\le R\log^4(K/\vare),
\end{equation}
and for all  $\CC \in \CK(\vare)$
\begin{equation}\label{5}
|(\BX_{\bsymb{\a}}\CC,\CC)_\CG|\le\vare\|\CC\|^2_\CG=\vare(|\CC^+\|^2_\CG+|\CC^-\|^2_\CG).
\end{equation}

Choose $\vare \in (0,1)$ to be such that $\a_{\pm}/(1-\vare) <
\sqrt2 \nu_{\pm}$, or equivalently,
\begin{equation}\label{8a}
    \y_{\pm} > \vare \mu_{\pm}.
\end{equation}
Let $\CL_{\pm}(\vare)$ be subspaces of $\CG_{\pm}$ of co-dimension
$N_-(0; \BI_\pm+\frac{\a_{\pm}}{1-\vare}\BB''_{\pm})$ which are
such that
\[
\|\CC^\pm\|^2_\CG+\frac{\a_{\pm}}{1-\vare}\bb''_{\pm}[\CC^{\pm}]\ge 0,
\qquad \forall\CC^{\pm} \in \CL_{\pm}(\vare).
\]
Then, for $\CC\in \left(\CL_+(\vare) \oplus \CL_-(\vare)\right)
\cap \CK(\vare)$ the inequality \eqref{5.2a} is satisfied. It follows
that
\begin{equation} \label{8}\begin{split}
&F(\y_+,\y_-) :=  N_-(0; \BI+\a_+\BB'_++\a_-\BB'_-) \\
 \le & N_-(0; \BI_++\frac{\a_+}{1-\vare}\BB''_+) + N_-(0;\BI_-+
\frac{\a_-}{1-\vare}\BB''_-) + R\log^4(K/\vare).
\end{split}\end{equation}
By (\ref{5.2}), this gives
\begin{equation}\label{9}\begin{split}
& F(\y_+,\y_-) \le  \left\{ M+\varphi_+((1-\vare)\mu_+)\right\}(\y_+-\vare \mu_+)^{-1/2}  \\
   &+\left\{ M+\varphi_-((1-\vare)\mu_-)\right\}(\y_--\vare
   \mu_-)^{-1/2}+R\log^4(K/\vare).\end{split}
\end{equation}
The inequalities \eqref{8a} guarantee that the estimate
\eqref{5.2} with $\mu_\pm$ replaced by $(1-\vare)\mu_\pm$ and,
correspondingly, $\y_\pm$ replaced by $\y_\pm -\vare \mu_\pm$ is
still valid.

Now we choose $\vare$, keeping in mind to optimize the right-hand
side in \eqref{9}. Let $\Psi(t)$ be a function described in
\eqref{add4}. Since $\psi(t) = o(t^{-1/4}) $, on choosing
\[
\vare = \vare(\mu_+,\mu_-)=\frac{1}{2} \Psi\left( \rm{min}\{\frac{\y_+}{\mu_+},
\frac{\y_-}{\mu_-}\}\right),
\]
we find that the inequalities (\ref{8a}) are satisfied. For if
$\y_+/\mu_+ \le \y_-/\mu_-$, then
\[
\vare < \Psi(\y_+/\mu_+) \le \Psi(\y_-/\mu_-) \le \y_+/\mu_+ \le
\y_-/\mu_-.
\]
Also, $\vare = o(\y_{\pm})$ as $\y_{\pm} \rightarrow 0$.
\vskip0.2cm
Introduce the function
\[\varphi(\mu_+,\mu_-)=\varphi_+((1-\vare)\mu_+)\frac{\y_+^{1/2}}{(\y_+ -\vare\mu_+)^{1/2}}
+\varphi_-((1-\vare)\mu_-)\frac{\y_-^{1/2}}{(\y_-
-\vare\mu_-)^{1/2}}.\] It is well-defined for
$(1-\mu_+^{-1},1-\mu_-^{-1})\in\Om_\Psi$ and
${\varphi}(\mu_+,\mu_-) \rightarrow 0 $ as $\mu_{\pm}\rightarrow
1$. The inequality \eqref{9} turns into
\begin{equation*}
 F(\y_+,\y_-)  \le  M \y_+^{-1/2} +M \y_-^{-1/2}
   +{\varphi}(\mu_+,\mu_-)(\y_+^{-1/2} + \y_-^{-1/2}) +
   R \log^4(K/\vare).
\end{equation*}
By \eqref{add4}, the last term here is $o(\y_{\pm}^{-1/2})$ and so
\begin{equation*}
F(\y_+,\y_-)   \le  M \y_+^{-1/2} +M \y_-^{-1/2}
   + \Phi(\mu_+,\mu_-)(\y_+^{-1/2} + \y_-^{-1/2})
\end{equation*}
where $ \Phi$ is a bounded function, defined on the same domain as $\varphi$ and
having the
same properties. The estimate is uniform for
$(\y_+/\mu_+,\y_-/\mu_-) \in \Omega_{\Psi}$.

To obtain the lower estimate we again choose $\CK(\vare)$ as in
(\ref{4}). There is a subspace $\CL(\vare)$ of $\CG$ of
co-dimension $N_-(0;\BI+\a_+\BB'_++\a_-\BB'_-)$ which is such that
\[
\|\CC\|^2_\CG+\a_+\bb'_+[\CC]+ \a_-\bb'_-[\CC]\ge 0,\qquad
\forall \CC\in \CK(\vare).
\]
Then, for $\CC\in \CG(\vare) \cap \CK(\vare)$,
\[
\|\CC\|^2_\CG+\a_+\bb''_+[\CC]+ \a_-\bb''_-[\CC]\ge\vare\|\CC\|^2_\CG.
\]
It follows that
\[\begin{split}
&N_-(0;\BI_++\frac{\a_+}{1+\vare}\BB''_+) + N_-(0;\BI_-
+\frac{\a_-}{1+\vare}\BB''_-)\\ &\le
N_-(0;\BI+\a_+\BB'_++\a_-\BB'_-) + R\log^4(K/\vare), \end{split}\]
and so
\begin{equation*}
    F(\y_+,\y_-) \ge N_-(0;\BI_++\frac{\a_+}{1+\vare}\BB''_+) + N_-(0;\BI_-
+\frac{\a_-}{1+\vare}\BB''_-)-R\log^4(K/\vare).
\end{equation*}
The rest of the argument is the same as for the upper estimate.
Actually it is easier, for if $\mu_{\pm}(\vare) := \sqrt2
\nu_{\pm}(1+\vare)/\a_{\pm}$, then
$(1-\mu_+^{-1}(\vare),1-\mu_-^{-1}(\vare)) $ automatically lies in
$\Omega_{\Psi}$ if $(1-\mu_+^{-1},1-\mu_-^{-1})$ does, and $
\mu_{\pm}(\vare) \rightarrow 1 $ as $\mu_{\pm} \rightarrow 1$.

All in all we have therefore shown that there exists a bounded
function $ \Phi(\mu_+,\mu_-)$ on $ \Om_{\Psi}$ which vanishes as $
(\mu_+,\mu_-) \rightarrow (1,1)$ and such that, uniformly for
$(\y_+/\mu_+,\y_-/\mu_-) \in \Omega_{\Psi}$,
\begin{equation*}
\left|N_-(r_{0,0}; \BA'_{\bsymb{\a},\bsymb{\nu}}) - M \y_+^{-1/2}
- M \y_-^{-1/2}\right| \le \Phi(\mu_+,\mu_-)\left( \y_+^{-1/2} +
\y_-^{-1/2}\right).
\end{equation*}
The proof of Theorem \ref{ESII1} is therefore complete.

\section{Acknowledgments}
The work on the paper started in the Summer of 2004 when one of
the authors (M.S.) was a guest of the School of Mathematics,
Cardiff University. M.S. takes this opportunity to express his
gratitude to the University for its hospitality and to the EPSRC
for financial support under grant GR/T01556.

\bibliographystyle{amsalpha}

\begin{thebibliography}{22}
\bibitem{BS} M.~Sh.~Birman and M.~Solomyak, {\it Schr\" odinger
operator. Estimates for number of bound states as function-theoretical
problem},
Spectral theory of operators (Novgorod, 1989), 1--54. English translation:
Amer. Math. Soc. Transl. Ser. 2, {\bf 150}, Amer. Math. Soc., Providence,
RI, 1992


\bibitem{EE} D.~E.~Edmunds and W.~D.~Evans, {\it Spectral Theory
and Differential Operators}, Oxford University Press, Oxford 1987.

\bibitem{ESI} W.~D.~Evans and M.~Solomyak, {\it Smilansky's model of irreversible
quantum graphs: I. The absloutely continuous spectrum,} Journal of
Physics, A: Mathematics and General. {\bf 38} (2005), 1-17.


\bibitem{GK} I.C. Gohberg and M.G. Krein, {\it Introduction to the theory
of linear non-selfadjoint operators in Hilbert space.} Izdat.
``Nauka'', Moscow 1965. English translation: Amer. Math.
Soc., Providence (1969).

\bibitem{Ku} P.~Kuchment, {\it Graph models for waves in thin structures,}
Waves Random Media  {\bf 12}  (2002),  no. 4, R1--R24.

\bibitem{SM} U.~Smilansky, {\it Irreversible quantum graphs}, Waves in
Random Media, {\bf 14} (2004), 143 -- 153.

\bibitem{S1} M.~Solomyak, {\it On a differential operator appearing
in the theory of irreversible quantum graphs}, {\it Waves in
Random Media}, {\bf 14} (2004), 173-185.

\bibitem {S2} M.~Solomyak, {\it On the discrete spectrum of a
family of differential operators}, Funct. Analysis and its appl.,
{\bf 38} (2004), 217-223.

\bibitem {S3} M.~Solomyak, {\it On a mathematical model of the
irreversible quantum graph}, St.-Petersburg Math. J., {\bf 17} (2005), in press.


\end{thebibliography}

\end{document}